\documentclass[11pt, a4paper]{article}
\usepackage{latexsym}
\usepackage{indentfirst}
\usepackage{graphicx}
\usepackage{amsmath}
\usepackage{amssymb}

\begin{document}
\title {\textbf{A modelling system for seaport activities}}
\author{Alina Costea, Ionela Rodica \c Ticu, Gheorghe Mishkoy }
\date{}
\maketitle

\begin{abstract}
In this paper we shall approach a modelling system for seaport activities based on the average waiting time and average queue length of ships in the seaport. We shall propose some suggestions for deepening and expanding this modelling system.\\

\end{abstract}

\section{Introduction}
A seaport system must be planned so as to satisfy the port facility o-perators, as well as the shipping companies. Ships must be received in the seaport at a minimum queuing time and a full capacity to use the facilities of the berth. As for the arrival times of ships and the receiving times of ships, the queuing theory is useful in analyzing and planning the seaport systems [1]. With the help of this theory, we can calculate the average waiting time of ships, the optimum number of berths in a seaport and the maximum number of ships that may wait in a seaport.	

\section{The average waiting time and average number of ships in the harbour system }

We shall consider a system with $A$ subsystems, each subsystem having $B_i$ seaports. The seaport $j$ of subsystem $i (j=1, 2,\ldots, B_i, i=1,2,\ldots ,A)$ has $S_{i,j}$ berths. We suppose that the arrival process of ships in the system is of Poisson type with parameter $\lambda$ and the time for taking over ships in berths is an independent random variable that has a probability distribution of Erlang type and is independent of the arrival processes.

	Let $u_{i,j,k} (k=1, 2, \ldots ,S_{i,j})$ be the average speed of receiving ships in the seaport in berth $k$ of seaport $j$ of the subsystem. Arrival processes are determined according to the terms of probability that have the following arrival times: $p_1\lambda, p_2\lambda,\ldots, p_A\lambda$. The $p_i$ probability can be determined as follows:

\begin{equation}\label{MCTec1}
p_i=\frac{\sum\limits_{j=1}^{B_i}S_{i,j}}{\sum\limits_{i=1}^{A}\sum\limits_{j=1}^{B_i}S_{i,j}}
\end{equation}

We suppose that the system is in a steady state and according to Yue Wuyi, Gu Jifa, Nie Jiayu, Hu Pingxian [3], the steady state solutions exist if $\rho<1$, where $\rho$ represents the  occupancy factor of the system and is given by:
\begin{equation}\label{MCTec2}
\rho=\frac{\lambda}{\sum\limits_{i=1}^{A}\sum\limits_{j=1}^{B_i}\sum\limits_{k=1}^{S_{i,j}}u_{i,j,k}}
\end{equation}

The arrival times in seaport $j$ of subsystem $i$ are forming a Poisson process with parameters $r_{i,j}$, where
\begin{equation}\label{MCTec3}
r_{i,j}=\frac{\lambda S_{i,j}}{\sum\limits_{j=1}^{B_i}S_{i,j}}
\end{equation}

Let $\rho_{i,j}$ be the occupancy factor of seaport $j$ of subsystem $i$, with the formula:
\begin{equation}\label{MCTec4}
\rho_{i,j}=\frac{r_{i,j}}{\sum\limits_{k=1}^{S_{i,j}}u_{i,j,k}}
\end{equation}

First, we find the formulas to calculate the average queue length and average waiting time of ships to subsystem $i$ and then we shall apply these formulas to the entire system.

A ship is allowed to enter seaport $j$ of subsystem $i$ and she finds at her arrival at least one of the berths $S_{i,j}$ free and gets immediately a free berth chosen randomly, so the ship is waiting in queue $k (k = 1,2,\ldots, R)$ of harbour $j$ until a berth is available. Ships will be received in berths under the rule: first input first output. A berth may not be unoccupied when ships are waiting in line. A ship leaves the system after it is completely taken over.

By applying the formulas suggested by Lee and Longton [2] and from Little's formula, we shall obtain the average waiting time and average queue length for ships that are in the system.

Let $E[W]_{i,j}$ be the average waiting time of ships in seaport $j$ of subsystem $i$ and let  $E[A]_{i,j}$ be the average number of ships waiting in seaport $j$ of subsystem $i$. The average waiting time $E[W]_{i,j}$ is obtained from the following formula:
$$E[W]_{i,j}=\frac{1+c^2}{2}E[W^*]_{i,j}$$
where
$$E[W^*]_{i,j}=\frac{(S_{i,j}\rho_{i,j})^{S_{i,j}}}{S_{i,j}!\sum\limits_{k=1}^{S_{i,j}}u_{i,j,k}(1-\rho_{i,j})^2}\left\{\sum\limits_{n=0}^{S_{i,j}-1}\frac{(S_{i,j}\rho_{i,j})^n}{n!}+\frac{(S_{i,j}\rho_{i,j})^{S_{i,j}}}{S_{i,j}!(1-\rho_{i,j})}\right\}^{-1}$$

By using Little's formula, we shall obtain the average number $E[n]_{i,j}$ of ships found in seaport $j$ of subsystem $i$:
$$E[n]_{i,j}=r_{i,j}E[W]_{i,j}$$

We note by $Var[G]$, the dispersion of random variation $G$, and \\$Var[G]=E[G^2]-(E[G])^2$ and by $c$, the coefficient of variation of the random variable $G$, $c=\dfrac{\sqrt{Var[G]}}{E[G]}=\sqrt{\dfrac{1}{n}}$ where $n$ is the degree of Erlang distribution.

The average number of ships in seaport $j$ of subsystem $i$ is:
$$E[Q]_{i,j}=E[n]_{i,j}+S_{i,j}\rho_{i,j}$$

The total average waiting time of the system is:
$$E[W]=\sum\limits_{i=1}^{A}\sum\limits_{j=1}^{B_{i}}E[W]_{i,j}$$

The total average time of ships waiting is:
$$E[n]=\sum\limits_{i=1}^{A}\sum\limits_{j=1}^{B_{i}}E[n]_{i,j}$$

The average number of ships in the system is: 
$$E[Q]=\sum\limits_{i=1}^{A}\sum\limits_{j=1}^{B_{i}}(E[n]_{i,j}+S_{i,j}\rho_{i,j})$$

\section{Conclusions and suggestions}
As we see, the above mentioned modelling system operates only in terms of average values and is valid only for modelling stationary characteristics. The non-stationary case, quite frequent in contemporary practice issues, including in port issues, remains outside the modelling system. According to restrictions p.2, the service time of ships is deemed to be distributed by the Erlang distribution, which also decreases the wide applicability of this system because the real service time may have some other distribution. A more efficient modelling system can be developed based on analytical results published in the monograph [4]. The models presented in this paper are invariant to the type of law for distributing the service, and characteristics obtained are valid for arbitrary distribution, which also offers opportunities for mode-lling the distribution time of service for a wide variety of distribution laws. Moreover, the analytical results are obtained for the non-stationary case, which offers us opportunities for modelling various extreme scenarios. The stationary case is obtained by passing to the limit and applying the Tauber theorem.

\end{document}